\documentclass[runningheads,a4paper]{llncs}

\usepackage{amssymb}
\usepackage{amsmath}
\usepackage{graphicx}

\usepackage{epsfig}
\usepackage{graphicx}
\usepackage{dcolumn}
\usepackage{bm}
\usepackage{tikz}
\usepackage[all]{xy}
\usetikzlibrary{arrows,shapes}
\usepackage{url}
\usepackage{centernot}
\usepackage{mathtools}
\usepackage{stmaryrd}

\usepackage{verbatimbox}

\begin{document}

\mainmatter  

\title{The multiresolution analysis of flow graphs
}

\titlerunning{The multiresolution analysis of flow graphs}

%
%
\author{Steve Huntsman\inst{1}
}
\authorrunning{S. Huntsman}

\institute{
BAE Systems FAST Labs, 4301 North Fairfax Drive, Arlington, VA 22203, USA\\
\email{steve.huntsman@baesystems.com}
}

%
%

\toctitle{}
\tocauthor{}
\maketitle

\begin{abstract}
We introduce and prove basic results about several graph-theoretic notions relevant to the multiresolution analysis of flow graphs that represent the transfer of control in computer programs. We take a category-theoretical viewpoint to demonstrate that our definitions are natural and to motivate particular incarnations of related constructions. 
\end{abstract}

\section{\label{sec:Introduction}Introduction}

The notion of a ``flow graph'' is central to the analysis and compilation of computer programs, encompassing constructs that represent the transfer of control and data \cite{CooperTorczon,Muchnick,NielsonNielsonHankin}. As the complexity of software increases, so does the scale of the corresponding flow graphs: accordingly, a framework for the analysis of flow graphs at multiple resolutions is desirable. Such a framework was originally presented in \cite{JohnsonPearsonPingali}, based on a hierarchical representation of input/output structure called the \emph{program structure tree} (PST). 

\begin{figure}
\begin{minipage}{0.3\textwidth}
\begin{verbnobox}[\scriptsize\arabic{VerbboxLineNo}\scriptsize\hspace{3ex}]
START
repeat
    repeat
        repeat
            if b goto 7
            if b
                repeat
                    S
                until b
            endif
        until b
        do while b
            do while b
                repeat
                    S
                until b
            enddo
        enddo
    until b
until b
HALT
\end{verbnobox}
\end{minipage}
\begin{minipage}{0.1\textwidth}
\ \\
\end{minipage}
\begin{minipage}{0.6\textwidth}
	\includegraphics[trim = 65mm 100mm 55mm 100mm, clip, width=\textwidth, keepaspectratio]{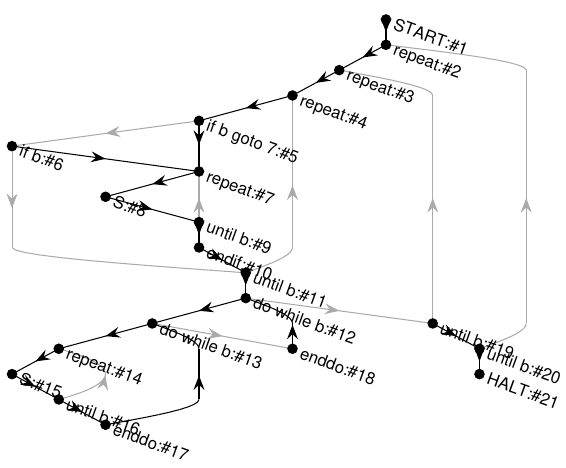}
\end{minipage}
	\caption{\label{fig:example1} (L) A simple imperative program. \texttt{S} denotes a generic statement (or subroutine); \texttt{b} denotes a generic Boolean predicate. (R) The corresponding control flow graph: branches are shaded {\color{black}black} (resp., {\color{gray}gray}) if the corresponding \texttt{b} evaluates to ${\color{black}\top}$ or ${\color{gray}\bot}$.}
\end{figure}

\begin{figure}
\begin{minipage}{0.5\textwidth}
\includegraphics[trim = 65mm 100mm 55mm 100mm, clip, width=\textwidth, keepaspectratio]{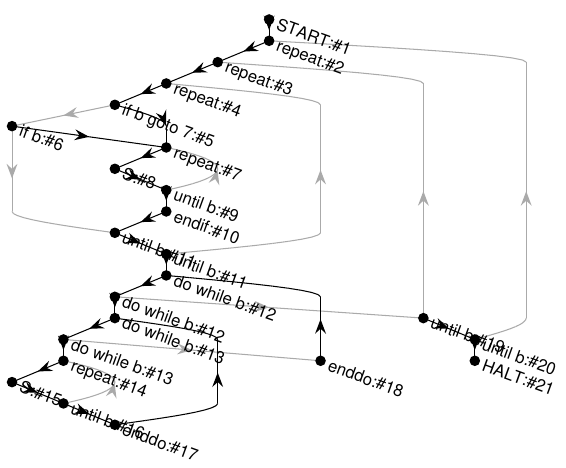}
\end{minipage}
\begin{minipage}{0.5\textwidth}
	\includegraphics[trim = 65mm 100mm 35mm 100mm, clip, width=\textwidth, keepaspectratio]{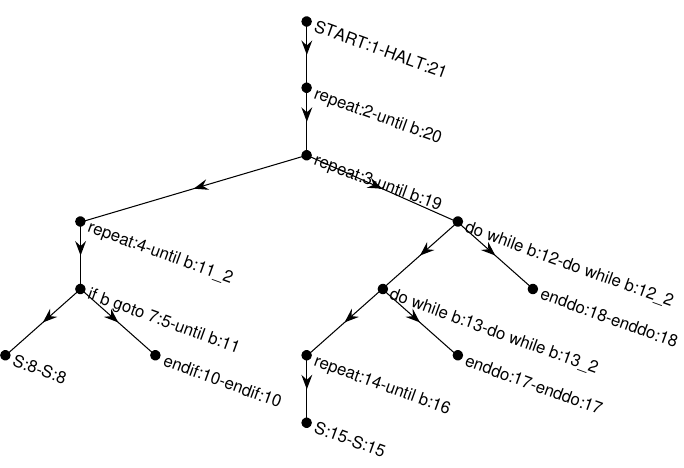}
\end{minipage}
	\caption{\label{fig:example1b} (L) ``Stretching'' (\emph{\`a la} \S \ref{sec:Stretching}) the flow graph of Figure \ref{fig:example1}. (R) The resulting PST.}
\end{figure}

\begin{figure}
\begin{centering}
	\includegraphics[trim = 60mm 95mm 55mm 100mm, clip, width=.3\textwidth, keepaspectratio]{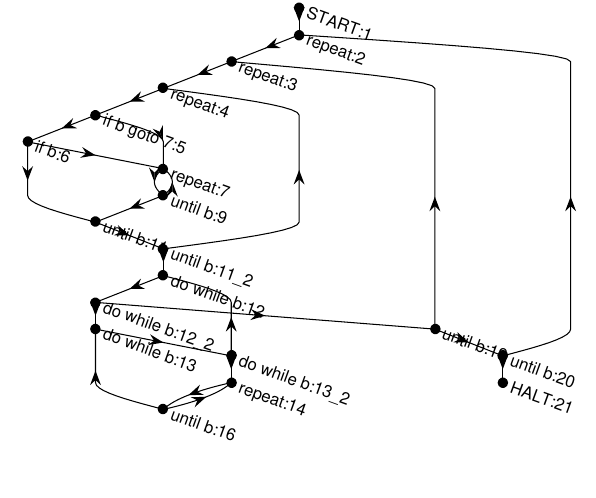}
	\includegraphics[trim = 60mm 95mm 55mm 100mm, clip, width=.3\textwidth, keepaspectratio]{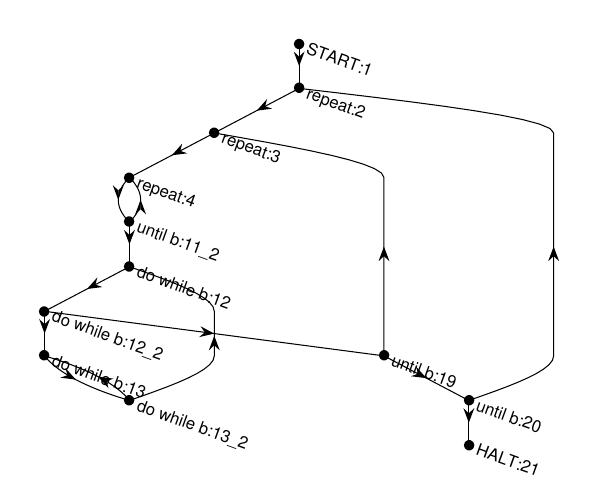}
	\includegraphics[trim = 60mm 95mm 55mm 100mm, clip, width=.3\textwidth, keepaspectratio]{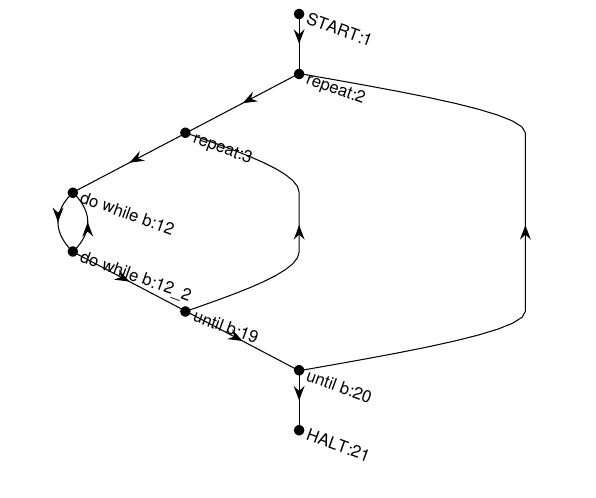}
	\caption{\label{fig:example1c} Succesively coarsening (\emph{\`a la} \S \ref{sec:Coarsening}) the flow graph of Figure \ref{fig:example1b}.}
\end{centering}
\end{figure}

For an illustration of this framework, consider the simple imperative program ``skeleton'' and associated control flow graph in Figure \ref{fig:example1}. The result of ``stretching'' it \emph{\`a la} \S \ref{sec:Stretching} and the PST of the result are shown in Figure \ref{fig:example1b}. Iterating the process of pruning each leaf of the PST \emph{\`a la} \S \ref{sec:Coarsening} leads to ``coarsened'' control flow graphs such as those in Figure \ref{fig:example1c}.

The utility of this framework is enhanced by \cite{ZhangDHollander}, which shows how to restructure the control flow graph of a program in such a way that subroutines can be identified as programs in their own right using the control flow graph alone. This feeds naturally into a ``multiresolution analysis'' of recursively composing (resp. decomposing) a program from (resp. into) subprograms in a way that can help with building, understanding, and modifying large programs.

This paper extends the work of \cite{JohnsonPearsonPingali} while correcting both an error of definition (for interiors of single-entry/single-exit regions) found in \cite{BoissinotBriskDarteRastello}, and another subtler error in the original proof of Theorem \ref{theorem:SESENesting}, by unifying and formalizing several natural concepts relevant to the decomposition and construction of flow graphs. This has several benefits: as the most basic example, we provide a definition of flow graph that is slightly different than its other usual variants but that is mathematically more natural and well-behaved. This in turn leads to a simpler analogue of the ``refined process structure tree'' of \cite{PolyvyanyyVanhataloVolzer,VanhataloVolzerKoehler} and natural category-theoretic constructions. These include multiresolution operations that  approximate and/or refine flow graphs at multiple scales, as well as series and parallel operations that respectively embody sequential execution and if/else constructs in control flow.

While most of the results of this paper are conceptually straightforward and many are at least latent in the literature, few of them have been simultaneously formulated explicitly and mathematically. Indeed, the practical motivation for this paper is simply to show that the ``right'' definition of a flow graph entails all the obvious desiderata, particularly for treating subroutines as programs in their own right. As \S \ref{sec:TTGs} highlights, the precise ability to compose flow graphs in category-theoretically nice ways is novel (though it is obvious that such a thing should be possible somehow): the unit object presents the principal difficulty, and much of our effort is focused on this issue for the case of parallel composition. This compositionality can inform the internal representation of graphical data structures and techniques for their manipulation in binary program analysis platforms such as \cite{BrumleyJagerThanassisSchwartz} and program synthesizers \cite{Gulwani} as well as compilers. 

In particular, constraining the notion of a valid control flow graph to the one considered in the paper could confer an advantage from the point of view of precompilation or reuse/modification: our results give a recipe for inserting and combining precompiled code in a convenient way. In a similar vein, we may want to understand a disassembled binary by synthesizing a similar or equivalent program. After restructuring the control flow graph along the lines of \cite{ZhangDHollander} and performing some straightforward normalizations (see \S \ref{sec:Stretching}), we could construct the PST and attempt program synthesis for each of the subroutines corresponding to a leaf node. In particular, we could generate inputs and observe outputs to each of these subroutines, so that program induction is a viable fallback at each point. Recursively going up the PST, we (attempt to) get such a globally synthesized program, and our results indicate precisely how synthesized/induced programs of intermediate scale can be maintained and reasoned over.

In other words, the constructions of the paper can inform tools that blur the lines between compilation and decompilation. In particular, the central results of \S \ref{sec:Coarsening} and \S \ref{sec:Tensor} contain the technical details necessary to have confidence that intermediate representations of programs can be (de)composed in a mathematically principled way, offering a firm foundation for future tools. Although superficial errors in \cite{JohnsonPearsonPingali} and hitherto unrecognized categorical structure in the PST have hindered its use,
\footnote{
To illustrate this point, we quote liberally from \cite{BoissinotBriskDarteRastello}: ``Unfortunately, we discovered an error in the aforementioned proof regarding SSI [static single information] form...we discovered that this mistake had been made in an earlier paper as well, and that other mistakes had been made in several papers that built on SSI form. The goal of this article, therefore, is to clear up the mistakes to the greatest possible extent...The key mistake was...made by Johnson \emph{et al.} [1994], who introduced a data structure called the program structure tree (PST), which attempted to represent the structure of a control flow graph hierarchically.''
}
we believe that tools based on it can and should be built.

The paper is organized as follows: we discuss dominance relations in \S \ref{sec:Dominance}; flow graphs, single-entry/single-exit regions, and the PST in \S \ref{sec:FlowSesePst}; we introduce the structure of a category on flow graphs in \S \ref{sec:Category} (this delay is to connect the paper to prior work most clearly); we discuss multiresolution transformations on flow graphs in \S \ref{sec:Coarsening}; and in \S \ref{sec:Tensor} we discuss series and parallel composition of flow graphs in the context of formal tensor product structures. 
\S \ref{sec:TTGs} discusses two-terminal graphs before our concluding remarks in \S \ref{sec:Conclusion}. \S \ref{sec:Proofs} contains proofs and \S \ref{sec:Stretching} sketches a ``stretching'' operation that enhances the applicability of our constructions.

We remark at the outset that all graphs (and related objects) are assumed finite throughout this paper. By convention, digraphs are allowed to have loops from a vertex to itself. Given a vertex $v$ in a digraph, let $d^+_0(v)$, $d^-_0(v)$, and $d^0(v)$ respectively denote the number of incoming edges excluding any loop, the number of outgoing edges excluding any loop, and the number ($\le 1$) of loops at $v$. A vertex $v$ is a \emph{source} iff $d^+_0(v) = 0$ and a \emph{target} iff $d^-_0(v) = 0$, i.e., loops have no bearing on these properties.

\section{\label{sec:Dominance}Dominance relations}

Let $G$ be a digraph and $j,k \in V(G)$. We say that $j$ \emph{dominates} $k$, written $j \text{ dom } k$, iff every path from a source $s$ in $G$ to $k$ passes through $j$ \cite{CooperTorczon,Muchnick}. Define $D_{jk} = 1$ if $j \text{ dom } k$ and $D_{jk} = 0$ otherwise. Similarly, let $D^\dagger := D(G^*)$, where $G^*$ is the reversal or adjoint of $G$ with adjacency matrix $A^*$ and corresponding dominance relation $\text{dom$^\dagger$}$. If $D^\dagger_{jk} = 1$, i.e., if $j \text{ dom$^\dagger$ } k$, write $k \text{ pdom } j$ and say that \emph{$k$ postdominates $j$}. Both the dominance and postdominance relations extend to edges. The following two lemmas are straightforward.

\begin{lemma}
\label{lemma:DominanceOrdering}
For distinct edges $\{e_j\}_{j=1}^3$ in a digraph $G$, if $e_1 \text{ dom } e_3$ and $e_2 \text{ dom } e_3$, then either $e_1 \text{ dom } e_2$ or $e_2 \text{ dom } e_1$. Similarly, if $e_1 \text{ pdom } e_2$ and $e_1 \text{ pdom } e_3$, then either $e_2 \text{ pdom } e_3$ or $e_3 \text{ pdom } e_2$. \qed
\end{lemma}

\begin{lemma}
\label{lemma:DominanceCycle}
If $e_1 \text{ dom } e_2$ and $e_1 \text{ pdom } e_2$ with $e_1 \ne e_2$, then a path from a source to a target that traverses $e_2$ contains a cycle of the form $(e_1,\dots,e_2,\dots,e_1)$. \qed
\end{lemma}

We use Lemma \ref{lemma:DominanceCycle} to fix a subtle (and minor) error in a proof of Theorem \ref{theorem:SESENesting} that was originally presented by \cite{JohnsonPearsonPingali}. This helps us to rescue the framework of \cite{JohnsonPearsonPingali} in its entirety from the problems raised by \cite{BoissinotBriskDarteRastello}.

\section{\label{sec:FlowSesePst}Flow graphs, single-entry/single-exit regions, and the program structure tree}

A \emph{flow graph} $G$ is a digraph with exactly one source and exactly one target, such that there is a unique (entry) edge from the source and a unique (exit) edge to the target, and such that identifying the source of the entry edge with the target of the exit edge yields a strongly connected digraph. (We do not require the entry and exit edges to be distinct, e.g., if $|V(G)| = 2$.) 
\footnote{NB. One sometimes sees variants of the definition and naming of this particular sort of concept, for the latter most typically as ``flowgraph'', ``flowchart'', or ``flow chart''. Some concepts with the same name are technically quite different but ``spiritually'' viewed in a similar context, as, e.g., in the work of Manin \cite{DelaneyMarcolli,Manin2009}.
}

\begin{figure}
	\centering
	\begin{tikzpicture}[scale = 0.5, every node/.style={inner sep=0,outer sep=0,scale=0.75},->,>=stealth',shorten >=1pt]
		\node (G) at (3,1.5) {$G$};
		\node [draw,circle,minimum size=5mm] (v01) at (-2.5,0) {1};
		\node [draw,circle,minimum size=5mm] (v02) at (2,0) {2};
		\node [draw,circle,minimum size=5mm] (v03) at (4,0) {3};
		\node [draw,circle,minimum size=5mm] (v04) at (6,0) {4};
		\node [draw,circle,minimum size=5mm] (v05) at (0,-2) {5};
		\node [draw,circle,minimum size=5mm] (v06) at (2,-2) {6};
		\node [draw,circle,minimum size=5mm] (v07) at (4,-2) {7};
		\node [draw,circle,minimum size=5mm] (v08) at (6,-2) {8};
		\node [draw,circle,minimum size=5mm] (v09) at (0,-4) {9};
		\node [draw,circle,minimum size=5mm] (v10) at (2,-4) {10};
		\node [draw,circle,minimum size=5mm] (v11) at (4,-4) {11};
		\node [draw,circle,minimum size=5mm] (v12) at (6,-4) {12};
		\node [draw,circle,minimum size=5mm] (v13) at (0,-6) {13};
		\node [draw,circle,minimum size=5mm] (v14) at (2,-6) {14};
		\node [draw,circle,minimum size=5mm] (v15) at (4,-6) {15};
		\node [draw,circle,minimum size=5mm] (v16) at (6,-6) {16};
		\node [draw,circle,minimum size=5mm] (v17) at (0,-8) {17};
		\node [draw,circle,minimum size=5mm] (v18) at (2,-8) {18};
		\node [draw,circle,minimum size=5mm] (v19) at (4,-8) {19};
		\node [draw,circle,minimum size=5mm] (v20) at (8.5,-8) {20};
		\foreach \from/\to in {
			v01/v02, v02/v03, v03/v04, v04/v08, v05/v02, v06/v09, v07/v06, v08/v07, v08/v12, v09/v05, v10/v14, 
			v11/v10, v11/v15, v12/v11, v13/v09, v14/v11, v14/v15, v15/v12, v15/v16, v16/v19, v17/v13, v18/v17, v19/v18, v19/v20}
			\draw (\from) to (\to);
		\draw (v07) to [out=60,in=120,looseness=8] (v07);
		\draw (v12) to [out=-30,in=30,looseness=8] (v12);
		\draw (v12) to [out=-75,in=75,looseness=1] (v16);
		\draw (v16) to [out=105,in=-105,looseness=1] (v12);
		\draw [lightgray] (3.5,-0.5) rectangle (4.5,0.5);	
		\draw [lightgray] (5.5,-0.5) rectangle (6.5,0.5);	
		\draw [lightgray] (-0.5,-2.5) rectangle (0.5,-1.5);	
		\draw [lightgray] (1.5,-2.5) rectangle (2.5,-1.5);	
		\draw [lightgray] (3.5,-2.5) rectangle (4.5,-0.8);	
		\draw [lightgray] (1.5,-4.5) rectangle (2.5,-3.5);	
		\draw [lightgray] (-0.5,-6.5) rectangle (0.5,-5.5);	
		\draw [lightgray] (-0.5,-7.5) rectangle (0.5,-8.5);	
		\draw [lightgray] (1.5,-7.5) rectangle (2.5,-8.5);	
		\draw [lightgray,style=dashed] (3.4,-0.6) rectangle (6.6, 0.6);	
		\path [draw,lightgray,style=dashed] (1.4,-2.6) -- (4.6,-2.6) -- (4.6,-0.7) -- (3.4,-0.7) -- (3.4,-1.4) -- (1.4,-1.4) -- cycle;	
		\path [draw,lightgray] (1.4,-3.4) -- (1.4,-6.6) -- (6.6,-6.6) -- (6.6,-4.6) -- (7.2,-4.6) -- (7.2,-3.4) -- cycle;	
		\draw [lightgray,style=dashed] (-0.6,-5.4) rectangle (0.6,-8.6);	
		\draw [lightgray,style=dashed] (-0.7,-7.3) rectangle (2.7,-8.7);	
		\path [draw,lightgray,style=dashed] (-0.8,-8.8) -- (2.8,-8.8) -- (2.8,-7.2) -- (0.8,-7.2) -- (0.8,-5.2) -- (-0.8,-5.2) -- cycle;	
		\draw [lightgray] (-1.5,1) rectangle (7.5,-9);	
		\node (G) at (0*2+13,1.5) {$\text{PST}(G)$};
		\node [draw,circle,minimum size=5mm] (v01) at (-1*2+13,0*2-4) {1};
		\node [draw,circle,minimum size=5mm] (v02) at (0*2+13,0*2-4) {2};
		\node [draw,circle,minimum size=5mm] (v03) at (-.5*2+13,.8660*2-4) {3};
		\node [draw,circle,minimum size=5mm] (v04) at (0*2+13,1*2-4) {4};
		\node [draw,circle,minimum size=5mm] (v05) at (.5*2+13,.8660*2-4) {5};
		\node [draw,circle,minimum size=5mm] (v06) at (.8660*2+13,.5*2-4) {6};
		\node [draw,circle,minimum size=5mm] (v07) at (1*2+13,0*2-4) {7};
		\node [draw,circle,minimum size=5mm] (v12) at (.8660*2+13,-.5*2-4) {12};
		\node [draw,circle,minimum size=5mm] (v13) at (.5*2+13,-.8660*2-4) {13};
		\node [draw,circle,minimum size=5mm] (v17) at (0*2+13,-1*2-4) {17};
		\node [draw,circle,minimum size=5mm] (v18) at (-.5*2+13,-.8660*2-4) {18};
		\node [draw,circle,minimum size=5mm] (v10) at (1.7321*2+13,-1*2-4) {10};
		\foreach \from/\to in {
			v01/v02, v02/v03, v02/v04, v02/v05, v02/v06, v02/v07, v02/v12, v02/v13, v02/v17, v02/v18, v12/v10}
			\draw (\from) to (\to);
	\end{tikzpicture}
	\caption{\label{fig:SESE} (L) SESE regions of the flow graph $G$ are outlined in gray: minimal (resp., non-minimal) SESE region outlines are solid (resp., dashed). Locally maximal but not minimal SESE regions are $((2,3),(4,8))$, $((8,7),(6,9))$, and $((19,18),(13,9))$. (R) The PST encodes the nesting of minimal SESE regions. Nodes are labeled by the target of the incoming edge (with a ``phantom'' edge from $-\infty$ to the source). The sets $\{3,4\}$, $\{6,7\}$, and $\{13,17,18\}$ correspond to locally maximal SESE regions that could sensibly be ``aggregated'' by identifying the respective vertices and omitting any resulting loops: however, a more mathematically natural variant of this construction is discussed in \S \ref{sec:Coarsening}.}
\end{figure}
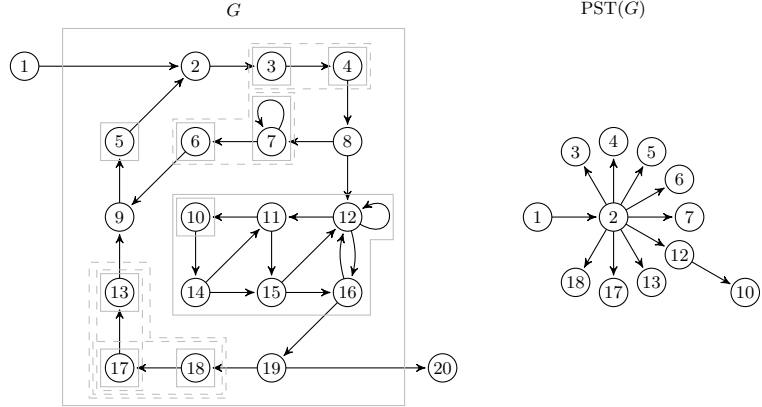

A \emph{single entry/single exit (SESE) region} in a digraph $G$ is defined as an ordered pair of edges $(e_1,e_2)$ satisfying each of the following conditions \cite{JohnsonPearsonPingali}: $e_1 \text{ dom } e_2$, $e_2 \text{ pdom } e_1$, and a cycle in $G$ contains $e_1$ iff it contains $e_2$.
See Figure \ref{fig:SESE} for examples. Note first that $(e_1,e_1)$ is a degenerate SESE region, 
\footnote{NB. Degenerate SESE regions $(e_1,e_1)$ are excluded by the original definition of \cite{JohnsonPearsonPingali}. We allow such regions to make the series tensor product of \S \ref{sec:TensorSeries} work nicely.
}
and second that a nondegenerate SESE region $(e_1,e_2)$ (i.e., a SESE region with $e_1 \ne e_2$) unambiguously corresponds to the ordered vertex pair $(t(e_1),s(e_2))$, where $s(\cdot)$ and $t(\cdot)$ respectively denote the source and target of an edge. We may use either the edge or vertex pairs above to specify a nondegenerate SESE region. Note also that in a DAG the third condition above is trivial. Finally, note that the edges $e_s$ from the source and $e_t$ to the target of a flow graph $G$ together define a SESE region and \emph{vice versa}. With this in mind, write either $G$ or $(e_s,e_t)$ for the flow graph or the equivalent SESE region. 

We give a few simple results (the first is straightforward enough that we omit a proof) before moving on to a fundamental theorem.

\begin{lemma}
\label{lemma:SESEConcatenation}
If $(e_1,e_2)$ and $(e_2,e_3)$ are SESE regions, then so is $(e_1,e_3)$. \qed
\end{lemma}

\begin{lemma}
\label{lemma:SESENesting}
If $(e_1,e_2)$ and $(e_1,e_3)$ are SESE regions with $e_2 \ne e_3$ and $e_2 \text{ dom } e_3$, then $(e_2,e_3)$ is a (nondegenerate) SESE region.
\end{lemma}


\begin{corollary} 
If $(e_1,e_2)$ is a SESE region with $e_2 \ne e_3$ and $e_2 \text{ dom } e_3$, and $(e_2,e_3)$ is not a SESE region, then $(e_1,e_3)$ is also not a SESE region. \footnote{A useful restatement of this is that if $(e_1,e_2)$ is a SESE region with $e_2 \ne e_3$ and $e_2 \text{ dom } e_3$, then $(e_1,e_3)$ is not a SESE region unless $(e_2,e_3)$ is a SESE region. 
} \qed
\end{corollary}

\begin{figure}
	\centering
	\begin{tikzpicture}[scale=1,every node/.style={inner sep=0,outer sep=0,scale=0.75},->,>=stealth',shorten >=1pt]
		\node [draw,circle,minimum size=5mm] (v01) at (0,0) {1};
		\node [draw,circle,minimum size=5mm] (v02) at (1,0) {2};
		\node [draw,circle,minimum size=5mm] (v03) at (2,0) {3};
		\node [draw,circle,minimum size=5mm] (v04) at (3,.25) {4};
		\node [draw,circle,minimum size=5mm] (v05) at (3,-.25) {5};
		\node [draw,circle,minimum size=5mm] (v06) at (4,0) {6};
		\node [draw,circle,minimum size=5mm] (v07) at (5,0) {7};
		\node [draw,circle,minimum size=5mm] (v08) at (6,.25) {8};
		\node [draw,circle,minimum size=5mm] (v09) at (6,-.25) {9};
		\node [draw,circle,minimum size=5mm] (v10) at (7,0) {10};
		\node [draw,circle,minimum size=5mm] (v11) at (8,0) {11};
		\node [draw,circle,minimum size=5mm] (v12) at (9,0) {12};
		\node [draw,circle,minimum size=5mm] (v13) at (4.5,.5) {13};
		\foreach \from/\to in {
			v01/v02, v02/v03, v03/v04, v03/v05, v04/v06, v05/v06, v06/v07, v07/v08, v07/v09, v08/v10, v09/v10, v10/v11, v11/v12}
			\draw (\from) to (\to);
		\draw (v11) to [out=150,in=0,looseness=1] (v13);
		\draw (v13) to [out=180,in=30,looseness=1] (v02);
	\end{tikzpicture}
	\caption{\label{fig:interiorFix} As \cite{BoissinotBriskDarteRastello} points out, the nondegenerate SESE regions $((2,3),(6,7))$ and $((6,7),(10,11))$ have interiors that intersect at vertex 13 according to the original definition of \cite{JohnsonPearsonPingali}. Our definition of the interior of a SESE region eliminates such unwanted behavior and allows us to salvage the original attempt to prove Theorem \ref{theorem:SESENesting}.}
\end{figure}
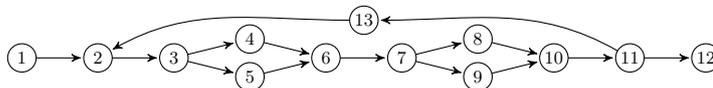

The \emph{interior} $G^\circ$ of $G \equiv (e_s,e_t)$ is the set of vertices that are each on at least one path starting from $t(e_s)$ that does not encounter $t(e_t)$. Critically, this definition differs slightly from Definition 6 of \cite{JohnsonPearsonPingali}, wherein the interior of a SESE region $(e_s,e_t)$ is defined as $\{j \in V : e_s \text{ dom } j \land e_t \text{ pdom } j\}$. An example in \S 5 of \cite{BoissinotBriskDarteRastello} and reproduced in Figure \ref{fig:interiorFix} illustrates the difference between these definitions.

A nondegenerate SESE region $(e_1,e_2)$ is called \emph{canonical} if for any SESE region $(e_1,e'_2)$ it is the case that $e_2 \text{ dom } e'_2$ and if for any SESE region $(e'_1,e_2)$ it is the case that $e_1 \text{ pdom } e'_1$. Our definition of the interior of a SESE region enables the following corrected version of Theorem 1 of \cite{JohnsonPearsonPingali} (cf. \cite{BoissinotBriskDarteRastello}).

\begin{theorem}
\label{theorem:SESENesting}
Interiors of distinct canonical SESE regions are disjoint or nested. 
\end{theorem}

Therefore canonical SESE regions are also \emph{minimal}, so we may use the two terms interchangeably: we generally prefer and use the latter. The inclusion relation on minimal SESE regions induces a tree---viz., the PST. An example of this nesting behavior and the corresponding PST are depicted in Figure \ref{fig:SESE}.

\begin{lemma}
\label{lemma:SESEDecomposition}
A nondegenerate SESE region $(e_0,e_\infty)$ decomposes as $(e_0,e_\infty) = \bigcup_{j = 1}^m (e_{j-1},e_j)$, where $e_m \equiv e_\infty$ and $(e_{j-1},e_j)$ are minimal SESE regions. 
\end{lemma}


Define an edge-indexed matrix $S$ by $S_{e_1,e_2} = 1$ if $(e_1,e_2)$ is a nondegenerate SESE region and $S_{e_1,e_2} = 0$ otherwise. Then $S$ is the adjacency matrix of a digraph whose weakly connected components correspond to the situation in Lemma \ref{lemma:SESEDecomposition}. We therefore obtain the following lemma.

\begin{lemma}
\label{lemma:Tournament}
Each weakly connected component of the digraph corresponding to $S$ is a transitive tournament, hence has a unique source, target, and a path of length 1 from source to target defining a \emph{locally maximal SESE region}. \qed
\end{lemma}

A closely related construction is the subject of \S \ref{sec:Coarsening}.

\section{\label{sec:Category}The category of flow graphs}

The principal goal of this section is merely to motivate and justify the details of the sequel. The key points are the introduction of the category ${\bf Dgph}$ of digraphs, and of its full subcategory ${\bf Flow}$ whose objects are flow graphs. 

It is natural to attempt to regard transformations of mathematical objects as morphisms in an appropriate category \cite{MacLane}. Unfortunately, in many if not most cases involving digraphs, such an attempt is complicated by technicalities that commonly arise from loops \cite{BrownMorrisShrimptonWensley}. The basic problem is that while identifying vertices should induce a graph morphism, such a morphism should also preserve edges. In particular, the morphism should preserve any edges between the vertices to be identified, necessarily inducing a loop. Insofar as we want loops in a coarse-grained control flow graph to correspond to actual loops in the atomic control flow, this is highly undesirable.

The common way around this problem is to treat loops on a separate footing. Following \cite{BrownMorrisShrimptonWensley}, define the category ${\bf Dgph}$ as follows. An object of ${\bf Dgph}$ is a \emph{reflexive digraph} $G = (U,\alpha,\omega)$ given by a set $U$ and \emph{head} and \emph{tail} functions $\alpha, \omega : U \rightarrow U$ satisfying $\alpha \circ \omega = \omega$ and $\omega \circ \alpha = \alpha$. Meanwhile, for $G' = (U',\alpha',\omega')$, a morphism $f \in {\bf Dgph}(G,G')$ is a function $f : U \rightarrow U'$ satisfying $f \circ \alpha = \alpha' \circ f$ and $f \circ \omega = \omega' \circ f$. 

The \emph{vertices} of $G = (U,\alpha,\omega)$ are the (mutual) image $V \equiv V(G)$ of $\alpha$ and $\omega$; the \emph{loops} are the set $L \equiv L(G) := \{u \in U : \alpha(u) = \omega(u)\}$ (so that $V \subseteq L$), and the \emph{edges} are the set $E \equiv E(G) := U \backslash L$.
\footnote{
The usual notion of a digraph is recovered by considering $\alpha \times \omega$ and its appropriate restrictions on $U^2$, $L^2$, and $E^2$: e.g., we can abusively write $E = (\alpha \times \omega)(E^2)$, where the LHS and RHS respectively refer to usual and reflexive notions of digraph edges. 
 } 
Thus a morphism $f : U \rightarrow U'$ restricts to $f|_V : V \rightarrow V'$, $f|_L : L \rightarrow L'$, and $f|_E : E \rightarrow U'$. 
In particular, morphisms are only partially specified by their actions on vertices, and the following definition is essentially a convention about how to treat vertex identification by default.

We define ${\bf Flow}$ to be the full subcategory of ${\bf Dgph}$ whose objects are (combinatorially realized as) flow graphs. 
\footnote{
As pointed out by D. Spivak, it would be desirable to describe flow graphs in terms of ${\bf Dgph}$, e.g. as algebras for some monad.
}

\section{\label{sec:Coarsening}Coarsening flow graphs}

We begin this section with intuition: the coarsening of a flow graph $G$ is obtained by taking each leaf of its PST and absorbing the interior of the corresponding sub-flow graph into its source. (See Figure \ref{fig:Coarsening}.) The details are below.

For $G \in {\bf Dgph}$, define the \emph{absorption} of $k$ into $j$ to be the morphism in ${\bf Dgph}$ (or the morphism's image, depending on context) which corresponds to identifying $k$ with $j$, and in the case $k \ne j$ subsequently annihilating any loop at $j$ (by mapping it to the vertex $j$). It is clear that first absorbing $k$ and then $m$ into $j$ is equivalent to first absorbing $m$ and then $k$ into $j$. Consequently, for $U \subseteq V(G)$ we may define the absorption of $U$ into $j$ in the obvious way. 
\footnote{
Failing to make fixed choices about whether to preserve or annihilate loops from, or formed at, absorbed and absorbing vertices amounts to a context-driven decision about the absorption process that is unlikely to be of any utility and need not be considered. Therefore, we proceed here to consider the space of such possible fixed choices. In the context of control flow graphs, a loop corresponds closely to a do-while construct. With this in mind, preserving such a construction under absorption corresponds to inserting additional computations into a do-while loop, or forming a new do-while loop around existing computations, altering the control flow. Meanwhile, annihilating loops corresponds to embedding the do-while construct within a larger sequence of computations, preserving the control flow. This is \emph{prima facie} cause to restrict consideration to the definition of absorption introduced above.
}

For $G, H \in {\bf Flow}$ with $H \subset G$, define the absorption of $H$ to be the result of absorbing the interior of $H$ into its source (considered as a vertex in $G$). This amounts to replacing $H$ with a single edge between its source and target. Finally, define the \emph{coarsening} $\circledcirc G$ of $G$ to be the result of absorbing all of the sub-flow graphs corresponding to leaves of the program structure tree of $G$. The fact that $\circledcirc G$ is well-defined follows from \cite{JohnsonPearsonPingali} (cf. the ``prime subprogram parse'' of \cite{TarjanValdes}) along with the preceding considerations. In particular, the definitions of absorption and coarsening yield the following technical lemma.

\begin{lemma}
\label{lemma:Coarsening}
Let $G \in {\bf Flow}$ and let $\circledcirc G$ result from absorbing the vertex sets $L_k$ into $k$ for all $k \in K$ (so that $L_k$ corresponds to a leaf of the program structure tree and $k \not \in L_k$). Let $L := \cup_{k \in K} L_k$ (this set should not be confused with the set of loops in $G$) and $J := V \backslash (K \cup L)$, so that $V = J \cup K \cup L$ and $J, K, L$ are mutually disjoint. Let $j, j' \in J$; $k, k' \in K$ with $k \ne k'$, and $\ell, \ell' \in L$. Finally, write $L_k^+ := \{k\} \cup L_k$ and let $g \in V$. Then the adjacency matrix of $\circledcirc G$ \emph{w.r.t. the vertex set of $G$} is $A'$, where $A'_{jj'} = A_{jj'}$, $A'_{jk'} = \bigvee_{\ell' \in L_{k'}^+} A_{j \ell'}$, $A'_{kj'} = \bigvee_{\ell \in L_k^+} A_{\ell j'}$, $A'_{kk'} = \bigvee_{\ell \in L_k^+, \ell' \in L_{k'}^+} A_{\ell \ell'}$, and $A'_{kk} = A'_{g \ell'} = A'_{\ell g'} = 0$. \qed
\end{lemma}

The real matter of substance in coarsening a flow graph is producing the sets $J$, $K$, and $L$ referred to just above (it turns out to be easier to construct the $L_k$ from $L$ than to go in the opposite direction).

\begin{theorem}
\label{theorem:Forest}
Using the notation of the preceding lemma, define a matrix $M$ as follows. For each leaf $(e_1,e_2)$ of the program structure tree, let $(e_1,e_2)^\circ$ denote its interior, and for all $j \in (e_1,e_2)^\circ$ set $M_{j,s(e_1)} = 1$. Then $M$ is the adjacency matrix of a DAG (in fact, a forest) whose weakly connected components have vertex sets $L_k^+$ and corresponding targets $k$.  
\end{theorem}


\begin{figure}
	\centering
	\begin{tikzpicture}[scale=0.5,every node/.style={inner sep=0,outer sep=0,scale=0.75},->,>=stealth',shorten >=1pt]
		\node (G) at (3,1) {$\circledcirc G$};
		\node [draw,circle,minimum size=5mm] (v01) at (0,0) {1};
		\node [draw,circle,minimum size=5mm] (v02) at (2,0) {2};
		\node [draw,lightgray,circle,minimum size=2mm] (v03) at (4,0) {};
		\node [draw,lightgray,circle,minimum size=2mm] (v04) at (6,0) {};
		\node [draw,lightgray,circle,minimum size=2mm] (v05) at (0,-2) {};
		\node [draw,lightgray,circle,minimum size=2mm] (v06) at (2,-2) {};
		\node [draw,lightgray,circle,minimum size=2mm] (v07) at (4,-2) {};
		\node [draw,circle,minimum size=5mm] (v08) at (6,-2) {8};
		\node [draw,circle,minimum size=5mm] (v09) at (0,-4) {9};
		\node [draw,lightgray,circle,minimum size=2mm] (v10) at (2,-4) {};
		\node [draw,circle,minimum size=5mm] (v11) at (4,-4) {11};
		\node [draw,circle,minimum size=5mm] (v12) at (6,-4) {12};
		\node [draw,lightgray,circle,minimum size=2mm] (v13) at (0,-6) {};
		\node [draw,circle,minimum size=5mm] (v14) at (2,-6) {14};
		\node [draw,circle,minimum size=5mm] (v15) at (4,-6) {15};
		\node [draw,circle,minimum size=5mm] (v16) at (6,-6) {16};
		\node [draw,lightgray,circle,minimum size=2mm] (v17) at (0,-8) {};
		\node [draw,lightgray,circle,minimum size=2mm] (v18) at (2,-8) {};
		\node [draw,circle,minimum size=5mm] (v19) at (4,-8) {19};
		\node [draw,circle,minimum size=5mm] (v20) at (6,-8) {20};
		\foreach \from/\to in {
			v01/v02, v02/v08, v08/v09, v08/v12, v09/v02, v11/v15, v12/v11, v14/v15, v15/v12, v15/v16, v16/v19, v19/v20}
			\draw (\from) to (\to);
		\draw (v19) to [out=150,in=-60,looseness=1.5] (v09);
		\draw (v12) to [out=-30,in=30,looseness=8] (v12);
		\draw (v11) to [out=-150,in=60,looseness=1] (v14);
		\draw (v14) to [out=30,in=-120,looseness=1] (v11);
		\draw (v12) to [out=-75,in=75,looseness=1] (v16);
		\draw (v16) to [out=105,in=-105,looseness=1] (v12);
		\draw [lightgray,style=dashed] (v03) to (v02);
		\draw [lightgray,style=dashed] (v04) to [out=-165,in=-15,looseness=1] (v02);
		\draw [lightgray,style=dashed] (v05) to (v09);
		\draw [lightgray,style=dashed] (v06) to [out=15,in=165,looseness=1] (v08);
		\draw [lightgray,style=dashed] (v07) to (v08);
		\draw [lightgray,style=dashed] (v10) to (v11);
		\draw [lightgray,style=dashed] (v13) to [out=-45,in=157.5,looseness=1] (v19);
		\draw [lightgray,style=dashed] (v17) to [out=15,in=165,looseness=1] (v19);
		\draw [lightgray,style=dashed] (v18) to (v19);
	\end{tikzpicture}
	\quad \quad \quad \quad
	\begin{tikzpicture}[scale=0.5,every node/.style={inner sep=0,outer sep=0,scale=0.75},->,>=stealth',shorten >=1pt]
		\node (G) at (3,1) {$\circledcirc^2 G$};
		\node [draw,circle,minimum size=5mm] (v01) at (0,0) {1};
		\node [draw,circle,minimum size=5mm] (v02) at (2,0) {2};
		\node [draw,circle,minimum size=5mm] (v08) at (6,-2) {8};
		\node [draw,circle,minimum size=5mm] (v09) at (0,-4) {9};
		\node [draw,lightgray,circle,minimum size=2mm] (v11) at (4,-4) {};
		\node [draw,lightgray,circle,minimum size=2mm] (v12) at (6,-4) {};
		\node [draw,lightgray,circle,minimum size=2mm] (v14) at (2,-6) {};
		\node [draw,lightgray,circle,minimum size=2mm] (v15) at (4,-6) {};
		\node [draw,lightgray,circle,minimum size=2mm] (v16) at (6,-6) {};
		\node [draw,circle,minimum size=5mm] (v19) at (4,-8) {19};
		\node [draw,circle,minimum size=5mm] (v20) at (6,-8) {20};
		\foreach \from/\to in {
			v01/v02, v02/v08, v08/v09, v08/v19, v09/v02, v19/v20}
			\draw (\from) to (\to);
		\draw (v19) to [out=150,in=-60,looseness=1.5] (v09);
		\draw [lightgray,style=dashed] (v11) to (v08);
		\draw [lightgray,style=dashed] (v12) to (v08);
		\draw [lightgray,style=dashed] (v14) to [out=60,in=-150,looseness=1] (v08);
		\draw [lightgray,style=dashed] (v15) to (v08);
		\draw [lightgray,style=dashed] (v16) to [out=75,in=-75,looseness=1] (v08);
	\end{tikzpicture}
	\caption{\label{fig:Coarsening} (L) Coarsening of the flow graph $G$ from Figure \ref{fig:SESE}. (Note that the pullback of the diagram $a \stackrel{g \circ f}{\longrightarrow} c \stackrel{g}{\longleftarrow} b$ is $a \stackrel{id}{\longleftarrow} a \stackrel{f}{\longrightarrow} b$, so that $f$ is the pullback of $g \circ f$ by $g$. We may therefore think of $\circledcirc G$ somewhat literally as a kind of pullback of $G$ by the leaves of its program structure tree.) (R) Coarsening again. A third coarsening is trivial.}
\end{figure}

Having considered coarsening flow graphs, we note that the appropriate mathematical formalization in the opposite direction---i.e., of inserting one flow graph into another \footnote{Note that we are not explicitly considering the insertion of loops in this setting.}---is captured by the assertion that flow graphs form a (symmetric) \emph{operad} \cite{Leinster,MarklShniderStasheff,Stasheff} (cf. \cite{RupelSpivak,Spivak2013}). 
At a high level, an operad is a collection of objects that ``plug into each other'' like maps $f_{(m)}: X^m \rightarrow X$ \emph{\`a la}
\begin{equation}
f_{(m)} \circ_\ell g_{(n)} := f(\cdot_1,\dots,\cdot_{\ell-1},g(\cdot_\ell,\dots,\cdot_{\ell+n-1}),\cdot_{\ell+n},\dots,\cdot_{m+n}). \nonumber
\end{equation} 
Let $P(n)$ denote the set of flow graphs with $n$ ordered edges and define the following family of maps
\begin{align}
\circ : P(n) \times P(k_1) \times \dots \times P(k_n) & \rightarrow P(k_1 + \dots + k_n) \nonumber \\
(G, G_1, \dots, G_n) & \mapsto G \circ (G_1, \dots, G_n)
\end{align}
by replacing, for each $1 \le j \le n$, the $j$th edge in $G$ with $G_j$ in the obvious way. Writing $k_0 \equiv 0$, the edge ordering on $G \circ (G_1, \dots, G_n)$ is obtained by assigning edges $\sum_{i=0}^{j-1} k_i + 1, \dots, \sum_{i=0}^j k_i$ to $G_j \hookrightarrow G \circ (G_1, \dots, G_n)$ in the same order as the edges of $G_j$, i.e., the edge ordering is inherited from its local components.

Definition-checking or direct comparison to other insertion operads (e.g. the little $d$-disks or $d$-cubes operads in ${\bf Top}$) yields the following 

\begin{theorem}
\label{theorem:operad}
The triple $\{e, \{P(n)\}_{n=1}^\infty, \circ\}$, where $e$ denotes the flow graph with one edge, forms an operad (in ${\bf Set}$). \qed
\end{theorem}

Thus the operadic composition $\circ$ and coarsening $\circledcirc$ operations are not only natural, but complementary, and we readily obtain the following lemma.

\begin{lemma}
\label{lemma:inverseCoarsening} 
If $G \in P(n)$ and $\circledcirc G_j = e \ne G_j$, then $\circledcirc (G \circ (G_1, \dots G_n)) = G$. \qed
\end{lemma}

\section{\label{sec:Tensor}Tensoring flow graphs}

\subsection{\label{sec:TensorSeries}Tensoring in series}

There is an essentially trivial tensor product on ${\bf Flow}$. The idea is simply to identify the exit edge of the first flow graph with the entry edge of the second flow graph, i.e., to combine flow graphs in series. The reason that this tensor product structure is interesting and useful is that it allows us a way to model additional structure in an enriched category. Specifically, this leads to the ${\bf Flow}$-category ${\bf SubFlow}_G$ of sub-flow graphs of a flow graph $G$.

We provide a quick sketch of the details here. Let $f \in {\bf Flow}(G,G_f)$ and $f' \in {\bf Flow}(G',G'_{f'})$ with $V(G) \cap V(G') = \varnothing$. Define $G \boxtimes G'$ to be the flow graph obtained by identifying the exit edge of $G$ and the entry edge of $G'$, and define $f \boxtimes f'$ to be the morphism in ${\bf Flow}(G \boxtimes G',G_f \boxtimes G'_{f'})$ obtained by identifying the output of $f$ on the exit edge of $G$ with that of $f'$ on the entry edge of $G'$. 

The following lemmas are straightforward. 

\begin{lemma}
\label{lemma:tensorSeries} 
${\bf Flow}$ is a monoidal category with tensor product given by $\boxtimes$, and with unit object the flow graph $e$ consisting of a single edge. \qed
\end{lemma}

\begin{lemma}
\label{lemma:enriched} 
For a generic flow graph $G$, we can form a category ${\bf SubFlow}_G$ enriched \cite{Kelly} over ${\bf Flow}$ as follows:
\begin{itemize}
\item ${\bf SubFlow}_G := E(G)$; \footnote{In particular, loops and reflexive self-edges are not included here, though the former may be accommodated without substantial changes.}
\item for $e_s, e_t \in {\bf SubFlow}_G$, the hom object ${\bf SubFlow}_G(e_s,e_t) \in {\bf Flow}$ is the (possibly empty) flow graph with entry edge $e_s$ and exit edge $e_t$; 
\item the composition morphism is induced by $\boxtimes$; 
\item the identity element is determined by the flow graph $e$ with one edge.  \qed
\end{itemize}
\end{lemma}

An important advantage of ${\bf SubFlow}_G$ over the path category of $G$ is that the former is finite (and the preceding sections essentially detail its construction), whereas the latter is infinite whenever there is a cycle in $G$.

\subsection{\label{sec:TensorParallel}Tensoring in parallel}

In this section we show that ${\bf Flow}$ carries a nontrivial monoidal structure (i.e., there is a tensor product operation that coherently combines flow graphs ``in parallel'' and not merely ``in series'' \cite{Coecke}). While the concept is rather obvious, the details are technical and we consequently make them explicit. In particular, although ${\bf Flow}$ is conceptually rather similar to the categories of $n$-cobordisms or tangles, the disjoint union only yields a tensor product in the latter cases: here, it must be modified to account for flow graphs whose entry and exit edges are identical or adjacent. 

Let $s(e^+), t(e^+), s(e^-), t(e^-)$ be four fixed distinct points not contained in the vertex set of any graph already under consideration, so that $e^\pm := (s(e^\pm),t(e^\pm))$ may be regarded as two separated abstract edges. If $G$ is a flow graph with entry edge $e_s$ and (possibly adjacent or even identical) exit edge $e_t$, define a ${\bf Dgph}$-morphism (i.e., the image may not be a flow graph) $\phi_G$ by the vertex/loop map
\begin{equation}
\label{eq:alphaG}
\phi_G(j) := 
\begin{cases}
   s(e^+) & \text{if } j = s(e_s) \\
   t(e^+) & \text{if } j = t(e_s) \\
   \ & \text{or } t(e_s) = s(e_t) \text{ and } j = t(e_t) \\
   s(e^-) & \text{if } t(e_s) \ne s(e_t) \text{ and } e_s \ne e_t \text{ and } j = s(e_t) \\
   t(e^-) & \text{if } t(e_s) \ne s(e_t) \text{ and } e_s \ne e_t \text{ and } j = t(e_t) \\
   j & \text{otherwise}
\end{cases}
\end{equation}
along with the extension to edges determined by not sending edges to (diagonal/reflexive edges or) loops. 

Intuitively, if the entry and exit edges of $G$ are neither identical nor adjacent, then $\phi_G$ maps them respectively to $e^+$ and $e^-$: otherwise, $\phi_G$ maps the entry edge to $e^+$ and everything else (vertices and edges to the vertex; loops to a loop) to $t(e^+)$. The rationale for the latter case is that it is the only really generic and consistent way for us to complete the definition of such a nontrivial ${\bf Dgph}$-morphism from a flow graph, and in fact this sort of definitional guidance is perhaps the primary rationale for invoking category theory \emph{ab initio}.

The following lemma is straightforward. 

\begin{lemma}
\label{lemma:phi} 
With  $j,k \in V(G)$ with $j \ne k$ and $j',k' \in V(G')$ with $j' \ne k'$, $\phi_G(j) = \phi_G(k) \Rightarrow \{j,k\} = \{t(e_s),t(e_t)\}$; similarly, $\phi_{G'}(j') = \phi_{G'}(k') \Rightarrow \{j',k'\} = \{t(e'_s),t(e'_t)\}$. \qed
\end{lemma}

If $V(G) \cap V(G') = \varnothing$, we define 
\begin{equation}
\label{eq:tensorObject}
G \otimes G' := G \sqcup G' / \sim, 
\end{equation}
where the equivalence relation on the disjoint (graph) union is determined, for $j,k \in V(G)$ with $j \ne k$ and $j',k' \in V(G')$ with $j' \ne k'$, by
\begin{align}
\label{eq:tensorSimilar}
(j,0) \sim (j,0) & \ \forall j; \nonumber \\
(j',1) \sim (j',1) & \ \forall j'; \nonumber \\
(j,0) \sim (k,0) & \iff (\phi_G(j) = \phi_G(k)) \wedge \star; \nonumber \\
(j',1) \sim (k',1) & \iff (\phi_{G'}(j') = \phi_{G'}(k')) \wedge \star'; \nonumber \\
(j,0) \sim (j',1) & \iff (\phi_G(j) = \phi_{G'}(j')) \wedge \star[j] \wedge \star[j'].
\end{align}
where we use the shorthands $\star := (t(e'_s) \ne s(e'_t)) \wedge (e'_s \ne e'_t)$; $\star' := (t(e_s) \ne s(e_t)) \wedge (e_s \ne e_t)$; $\star[j] := (t(e_s) = s(e_t)) \wedge (j \in \{t(e_s),t(e_t)\}) \Rightarrow \star$, and $\star[j'] := (t(e'_s) = s(e'_t)) \wedge (j' \in \{t(e'_s),t(e'_t)\}) \Rightarrow \star'$,
with an obvious extension to edges. (Here $e'_s$ and $e'_t$ denote the entry and exit edges of $G'$.)

\begin{lemma}
\label{lemma:tensorSimilar} 
\eqref{eq:tensorSimilar} indeed defines an equivalence relation.
\end{lemma}

The following lemma is straightforward. 

\begin{lemma}
\label{lemma:parallelTensorFlow} 
If $G$ and $G'$ are flow graphs, then so is $G \otimes G'$. Furthermore, if $e$ denotes the flow graph with a single edge, then $G \otimes e \cong e \otimes G \cong G$. \qed
\end{lemma}

Thus in particular we have inclusions $i_G : \phi_G(G) \hookrightarrow G \otimes G'$ and $i'_{G'} : \phi_{G'}(G') \hookrightarrow G \otimes G'$ given respectively by $i_G(\phi_G(j)) = [(j,0)]$ and $i'_{G'}(\phi_{G'}(j')) = [(j',1)]$, where as per usual practice $[\cdot]$ indicates an equivalence class under $\sim$. $G \otimes G'$ is a flow graph formed by identifying the entry edges of $G$ and $G'$, and identifying the respective exit edges if this does not affect the interiors of the factors, and otherwise collapsing the smaller factor in a way that sufficiently extends the identification of entry edges.

Meanwhile, for $f \in {\bf Flow}(G,G_f)$ and $f' \in {\bf Flow}(G',G'_{f'})$, we define $f \otimes f' \in {\bf Flow}(G \otimes G', G_f \otimes G'_{f'})$ as follows (see also Figure \ref{fig:tensorMapDiagram}): 
\begin{equation}
\label{eq:tensorMorphism}
(f \otimes f')(k) := 
\begin{cases}
   [(f(j),0)] & \text{if } k = [(j,0)] \\
   [(f'(j'),1)] & \text{if } k = [(j',1)]
\end{cases}
\end{equation}
along with the implied extension to edges.

\begin{figure}
	\centering
	\begin{tikzpicture}[->,>=stealth',shorten >=1pt]
		\node (v01) at (-5,0) {$G$};
		\node (v02) at (-2.5,0) {$\phi_G(G)$};
		\node (v03) at (0,0) {$G \otimes G'$};
		\node (v04) at (2.5,0) {$\phi_{G'}(G')$};
		\node (v05) at (5,0) {$G'$};
		\node (v06) at (-5,-1.5) {$G_f$};
		\node (v07) at (-2.5,-1.5) {$\phi_{G_f}(G_f)$};
		\node (v08) at (0,-1.5) {$G_f \otimes G'_{f'}$};
		\node (v09) at (2.5,-1.5) {$\phi_{G'_{f'}}(G'_{f'})$};
		\node (v10) at (5,-1.5) {$G'_{f'}$};
		\path [->] (v01) edge node [above] {$\phi_G$} (v02);
		\path [right hook->] (v02) edge node [above] {$i_G$} (v03);
		\path [left hook->] (v04) edge node [above] {$i'_{G'}$} (v03);
		\path [->] (v05) edge node [above] {$\phi_{G'}$} (v04);
		\path [->] (v01) edge node [left] {$f$} (v06);
		\path [->] (v03) edge node [left] {$f \otimes f'$} (v08);
		\path [->] (v05) edge node [right] {$f'$} (v10);
		\path [->] (v06) edge node [above] {$\phi_{G_f}$} (v07);
		\path [right hook->] (v07) edge node [above] {$i_{G_f}$} (v08);
		\path [left hook->] (v09) edge node [above] {$i'_{G'_{f'}}$} (v08);
		\path [->] (v10) edge node [above] {$\phi_{G'_{f'}}$} (v09);
	\end{tikzpicture}
	\caption{\label{fig:tensorMapDiagram} The tensor product of morphisms in ${\bf Flow}$.}
\end{figure}

\begin{lemma}
\label{lemma:tensorMorphismWellDefined}
\eqref{eq:tensorMorphism} is well-defined.
\end{lemma}


\begin{theorem}
\label{theorem:parallelTensor}
${\bf Flow}$ is a monoidal category with tensor product $\otimes$ given by \eqref{eq:tensorObject} and \eqref{eq:tensorMorphism}, and with unit object the flow graph $e$ consisting of a single edge.
\end{theorem}

\begin{corollary}
\label{theorem:parallelSMC}
$({\bf Flow},\otimes)$ is a symmetric monoidal category. \qed
\end{corollary}

\subsection{\label{sec:TensorRemark}Remarks}

The series and parallel tensor operations described above are very similar in spirit to the composition operations encountered in the study of so-called \emph{series-parallel graphs} \cite{BangJensenGutin,Duffin} (cf. \cite{DoughertyGutierrez}). While the category-theoretical analysis of series and parallel tensor operations in the context of something like a system or wiring diagram has a very long history \cite{Bainbridge}, a precise treatment appropriate to our development does not appear to be present in the literature.

\section{\label{sec:TTGs}Two-terminal graphs}

Many of the considerations of the present paper have obvious analogues in the case of \emph{two-terminal graphs} (TTGs). In particular, \cite{TarjanValdes,VanhataloVolzerKoehler} describes a multiresolution decomposition of TTGs (cf. \cite{Fujishige}) that is a more granular version of the PST. This \emph{refined process structure tree} reduces via a straightforward graph transformation (similar to that in \S \ref{sec:Stretching}) 
to constructing the SPQR tree \cite{PolyvyanyyVanhataloVolzer}. 

Unfortunately, computing SPQR trees is a notoriously intricate exercise: indeed, a correct linear time algorithm was not actually implemented until 2000 \cite{GutwengerMutzel}, though an incorrect version of the same algorithm was first described in 1973 \cite{HopcroftTarjan}. Today there is still not a completely explicit description of the correct linear time algorithm in the literature: for such an account it is necessary to refer to one of the two known publically available implementations in the C++ OGDF
\footnote{\url{http://www.ogdf.net/}}
and the Java jBPT
\footnote{\url{https://code.google.com/archive/p/jbpt/}} 
frameworks. 
\footnote{An alternative algorithm is in \cite{Tsin}, but we are not aware of an implementation.
} \footnote{For acyclic TTGs it is not hard to see that the analogue of SESE regions are vertex pairs $(j,k)$ s.t. $D_{jk} D^\dagger_{jk} - \delta_{jk} = 1$, but the cyclic case is much harder.
}

While the computation and properties of the fundamental decomposition for TTGs are more involved than the PST, some analogues of the constructions detailed in this paper are simpler since TTGs are defined to omit loops. On the other hand, one minor complication relative to flow graphs that informs notions of coarsening and inclusion operads for TTGs is that some TTGs can have their sources and targets swapped. A more significant (and perhaps surprising) complication is that it is not clear how to define a canonical parallel tensor operation for TTGs: the principal difficulty is the unit object. Lacking such an operation would be a significant shortcoming relative to the framework for flow graphs, as parallel tensoring in ${\bf Flow}$ corresponds to introducing an if/else statement in control flow.

\section{\label{sec:Conclusion}Conclusion}

Besides applications to understanding and manipulating programs mentioned in \S \ref{sec:Introduction}, our particular notion of a flow graph naturally yields an interesting category that readily admits explicit representations and manipulations in (and of!) software. While some of the constructions involved are somewhat delicate and inelegant (for example, much of \S \ref{sec:TensorParallel}), this is due to properly accounting for degenerate cases that are of little practical concern but that nevertheless constrain practical and principled techniques for representing, reasoning about, and composing program artifacts. 

Put another way, requiring that flow graphs exhibit category-theoretical desiderata places strong but satisfiable restrictions on them that can usefully inform the architecture of program analysis platforms, program synthesizers, compilers, etc. More generally, category theory allows us to address corner cases in the construction and manipulation of data structures whose resolution is not obvious.

\section*{\label{sec:acknowledgements}Acknowledgements}

We thank Brendan Fong, Artem Polyvyanyy, and David Spivak for helpful comments.

\appendix

\section{\label{sec:Proofs}Proofs}

\begin{proof} [Lemma \ref{lemma:SESENesting}]
The only thing to show is that $e_3 \text{ pdom } e_2$. It must be the case that either $e_2 \text{ pdom } e_3$ or $e_3 \text{ pdom } e_2$, so assume the former. Since $e_2 \text{ dom } e_3$ also, we must have that any source-target path traversing $e_3$ contains a cycle of the form $(e_2,\dots,e_3,\dots,e_2)$ by Lemma \ref{lemma:DominanceCycle}; deleting all cycles from this path yields a source-target path traversing $e_2$ but not $e_3$. Reversing this path yields a contradiction to the assumption that $e_2 \text{ pdom } e_3$. \qed
\end{proof}

\begin{proof} [Theorem \ref{theorem:SESENesting}]
[Although our definition of the interior of a SESE differs in a slight but critical way from from \cite{JohnsonPearsonPingali}, the proof is a mostly straightforward adaptation of the original attempt. That said, we also fix a minor gap of the original attempt for case ii).]

Let $(e_1,e_2)$ and $(e'_1,e'_2)$ be distinct canonical SESE regions whose interiors are not disjoint, and let $v$ be in their intersection. Since $e_1 \text{ dom } v$ and $e'_1 \text{ dom } v$, it must be that either $e_1 \text{ dom } e'_1$ or $e'_1 \text{ dom } e_1$: assume the former w.l.o.g. Similarly, since $e_2 \text{ pdom } v$ and $e'_2 \text{ pdom } v$, either $e_2 \text{ pdom } e'_2$ or $e'_2 \text{ pdom } e_2$: in the former case, $(e'_1,e'_2) \subset (e_1,e_2)$ and we are done, so assume the latter case. We now have three cases to consider: i) $e_2 = e'_1$; ii) $e_2 \ne e'_1$ and $e'_1 \text{ dom } e_2$; and iii) $e_2 \ne e'_1$ and $e'_1$ does not dominate $e_2$. We shall show that each case leads to a contradiction.

\emph{Case i).} Since in this case $e_2 = e'_1$, we have that $e_2 \text{ dom } v$ and $e_2 \text{ pdom } v$, so it must be that any path from the source to the target that traverses $v$ must contain a cycle of the form $(e_2,\dots,v,\dots,e_2)$ by Lemma \ref{lemma:DominanceCycle}. But this means that $v$ cannot be in the interior of $(e_1,e_2)$, a contradiction: hence case i) cannot hold.

\emph{Case ii).} Since in this case $e'_1 \text{ dom } e_2$ and generically $e_1 \text{ dom } e'_1$, we may decompose any path $\gamma_{02}$ from the source to $e_2$ (using an obvious notation) as $\gamma_{02} \equiv \gamma_{01} \gamma_{11'} \gamma_{1'2}$. Meanwhile since $e_2 \text{ pdom } e_1$, we may decompose any path $\gamma_{1 \infty}$ from $e_1$ to the target as $\gamma_{1 \infty} \equiv \gamma_{12} \gamma_{2 \infty}$. Taken together, these decompositions imply that we can decompose any path from the source to the target that traverses $e'_1$ as $\gamma_{01} \gamma_{11'} \gamma_{1'2} \gamma_{2 \infty}$, so that $e_2 \text{ pdom } e'_1$ and $e'_1 \text{ pdom } e_1$. 

Moreover, if there is a cycle that traverses $e_1$, it also traverses $e_2$ and \emph{vice versa}, so we may write such a cycle as $\omega_{12} \equiv \gamma_{12} \gamma_{21}$, where $\gamma_{12} \equiv \gamma_{11'} \gamma_{1'2}$ as above. Hence such a cycle $\omega_{12}$ must traverse $e'_1$. Similarly, if there is a cycle that traverses $e'_1$, it also traverses $e'_2$ and \emph{vice versa}, so we may write such a cycle as $\omega_{1'2'} \equiv \gamma_{1'2'} \gamma_{2'1'}$, where $\gamma_{1'2'}$ traverses $e_2$ since $e'_2 \text{ pdom } e_2$. Hence such a cycle $\omega_{1'2'}$ must traverse $e_2$. It follows that $(e'_1,e_2)$ is a SESE region. 

Since both $(e_1,e_2)$ and $(e'_1,e'_2)$ are canonical SESE regions, we have that $e_1 \text{ pdom } e'_1$ and $e'_2 \text{ dom } e_2$. At the same time, $e'_1 \text{ pdom } e_1$, so it must be that $e_1 = e'_1$. It follows that $(e_1,e'_2)$ is also a SESE region, and therefore also that $e_2 \text{ dom } e'_2$, so it must be that $e_2 = e'_2$. This contradicts the hypothesis that $(e_1,e_2)$ and $(e'_1,e'_2)$ are distinct: hence case ii) cannot hold.

\emph{Case iii).} Since in this case $e'_1$ does not dominate $e_2$, there is a path $\gamma_{02}$ from the source to $e_2$ that avoids $e'_1$. Suppose that $e'_1$ does not postdominate $e_2$, i.e., suppose that there is a path $\gamma_{2\infty}$ from $e_2$ to the target that avoids $e'_1$. Then since $e'_2 \text{ pdom } e_2$, $\gamma_{2\infty}$ must traverse $e'_2$. But since $e'_1 \text{ dom } e'_2$ and the concatenated path $\gamma \equiv \gamma_{02} \gamma_{2 \infty}$ from the source to the target traverses $e'_2$, it must be that $\gamma_{2\infty}$ traverses $e'_1$, contradicting the assumption that 
$e'_1$ does not postdominate $e_2$. Therefore since $e'_1 \text{ pdom } e_2$ and $e_2 \text{ pdom } v$, we have that $e'_1 \text{ pdom } v$. Moreover, $e'_1 \text{ dom } v$, so any path from the source to the target that traverses $v$ must contain a cycle of the form $(e'_1,\dots,v,\dots,e'_1)$ by Lemma \ref{lemma:DominanceCycle}. But this means that $v$ cannot be in the interior of $(e'_1,e'_2)$. By contradiction, case iii) cannot hold.
\qed
\end{proof}

\begin{proof} [Lemma \ref{lemma:SESEDecomposition}]
Suppose w.l.o.g. that $(e_0,e_\infty)$ is not minimal. Then at least one of the following is true: i) there exists a nondegenerate SESE region $(e_0,e_1)$ such that $e_\infty$ does not dominate $e_1$; ii) there exists a nondegenerate SESE region $(e_{-1},e_\infty)$ such that $e_0$ does not postdominate $e_{-1}$. Consider case i), and assume w.l.o.g. that $(e_0,e_1)$ is minimal (otherwise, we have at least one of case i) or ii) again). Then $e_1 \text{ dom } e_\infty$, so $(e_1,e_\infty)$ is a nondegenerate SESE region and we can write $(e_0,e_\infty) = (e_0,e_1) \cup (e_1,e_\infty)$. Exactly similar reasoning informs case ii), and an induction establishes the lemma. \qed
\end{proof}

\begin{proof} [Theorem \ref{theorem:Forest}]
 Let $(e_1,e_2)$ be a leaf of the PST. If $s(e_1)$ is in the interior of some other leaf $(e_1',e_2')$ of the PST, then $e_1 = e_2'$. Therefore, $M_{s(e_2'),s(e_1')} = 1$ and any other vertices $j$ with $M_{j,s(e_1')} = 1$ correspond to the remaining elements of $(e_1',e_2')^\circ$, which are leaves in the digraph $G_M$ with adjacency matrix $M$. On the other hand, if $s(e_1)$ is not in the interior of some other leaf of the PST, then it is a target in $G_M$. The result follows.  \qed
\end{proof}

\begin{proof} [Lemma \ref{lemma:tensorSimilar}]
Since it is obvious from the structure of $\star[j]$ and $\star[j']$ that $(j',1) \sim (j,0) \iff (j,0) \sim (j',1)$, the only thing to show is transitivity. A (perhaps unnecessarily) mechanical proof consists of verifying each of the eight assertions $(\ell_{1,b_1},b_1) \sim (\ell_{2,b_2},b_2) \sim (\ell_{3,b_3},b_3) \Rightarrow (\ell_{1,b_1},b_1) \sim (\ell_{3,b_3},b_3)$ for $(b_1,b_2,b_3) \in \{0,1\}^3$ and $\ell_{1,b_1}, \ell_{2,b_2}, \ell_{3,b_3}$ distinct. 

First, consider $(b_1,b_2,b_3) = (0,0,0)$: we must show in this case that $(\phi_G(\ell_{10}) = \phi_G(\ell_{20}) = \phi_G(\ell_{30})) \wedge \star$ implies
$(\phi_G(\ell_{10}) = \phi_G(\ell_{30})) \wedge \star$, 
but this is trivial.

Next, consider $(b_1,b_2,b_3) = (0,0,1)$. Here we must show that $(\phi_G(\ell_{10}) = \phi_G(\ell_{20}) = \phi_{G'}(\ell_{31})) \wedge \star \wedge \star[\ell_{20}] \wedge \star[\ell_{31}]$ implies $(\phi_G(\ell_{10}) = \phi_{G'}(\ell_{31})) \wedge \star[\ell_{10}] \wedge \star[\ell_{31}]$.
By Lemma \ref{lemma:phi}, $\{\ell_{10},\ell_{20}\} = \{t(e_s),t(e_t)\}$, so $t(e_s) = s(e_t)$ and $\star[\ell_{10}]$ is true, establishing the desired result.

For $(b_1,b_2,b_3) = (0,1,0)$, we must show that $(\phi_G(\ell_{10}) = \phi_{G'}(\ell_{21}) = \phi_G(\ell_{30})) \wedge \star[\ell_{10}] \wedge \star[\ell_{21}] \wedge \star[\ell_{30}]$
implies
$(\phi_G(\ell_{10}) = \phi_G(\ell_{30})) \wedge \star$. 
By Lemma \ref{lemma:phi}, $\{\ell_{10},\ell_{30}\} = \{t(e_s),t(e_t)\}$, so $t(e_s) = s(e_t)$ and $\ell_{10},\ell_{30} \in \{t(e_s),t(e_t)\}$. Since in the present case both $\star[\ell_{10}]$ and $\star[\ell_{30}]$ are true by assumption and we have just shown their hypotheses true, their mutual conclusion $\star$ is also true here. This yields the desired implication. (NB. Although $\star[\ell_{21}]$ is true in this case, neither its hypothesis nor its conclusion are.)

By symmetry, the last case we need to consider is $(b_1,b_2,b_3) = (0,1,1)$: we need to show here that
$(\phi_G(\ell_{10}) = \phi_{G'}(\ell_{21}) = \phi_{G'}(\ell_{31})) \wedge \star[\ell_{10}] \wedge \star[\ell_{21}] \wedge \star'$
implies
$(\phi_G(\ell_{10}) = \phi_{G'}(\ell_{31})) \wedge \star[\ell_{10}] \wedge \star[\ell_{31}]$. 
By Lemma \ref{lemma:phi}, $\{\ell_{21},\ell_{31}\} = \{t(e'_s),t(e'_t)\}$, so $t(e'_s) = s(e'_t)$ and $\star[\ell_{31}]$ is true, so we are done. \qed
\end{proof}

\begin{proof} [Lemma \ref{lemma:tensorMorphismWellDefined}]
We need to show that whenever $[(j,0)] = [(j',1)]$ we also have $[(f(j),0)] = [(f'(j'),1)]$. An equivalent assertion is that whenever $\phi_G(j) = \phi_{G'}(j')$, we also have $\phi_{G_f}(f(j)) = \phi_{G'_{f'}}(f'(j'))$. There are precisely four cases in which the hypothesis can hold, corresponding to the first four cases of \eqref{eq:alphaG} (note that the second case has four subcases). In the first case, both $[(f(j),0)]$ and $[(f'(j'),1)]$ must be the source of the entry edge in $G_f \otimes G'_{f'}$ since $f$ and $f'$ are morphisms in ${\bf Flow}$; similarly, the other cases respectively give that both $[(f(j),0)]$ and $[(f'(j'),1)]$ must be the target of the entry edge, the source of the exit edge, and the target of the exit edge. \qed
\end{proof}

\begin{proof} [Theorem \ref{theorem:parallelTensor}]
We must establish two things: that $\otimes$ is a bifunctor, and that it satisfies the necessary coherence conditions. 

To see that $\otimes$ is a bifunctor, first note that $(id_G \otimes id_{G'})([(j,0)]) = [(j,0)] = id_{G \otimes G'}([(j,0)])$
by \eqref{eq:tensorMorphism}, and 
$(id_G \otimes id_{G'})([(j',1)]) = [(j',1)] = id_{G \otimes G'}([(j',1)])$,
so that $id_G \otimes id_{G'} = id_{G \otimes G'}$. Now we must show that $(g \otimes g') \circ (f \otimes f') = (g \circ f) \otimes (g' \circ f')$. But this is easily seen since, again by \eqref{eq:tensorMorphism}, we have $(g \otimes g')([(f(j),0)]) = [(g(f(j)),0)] = ((g \circ f) \otimes (g' \circ f'))([(j,0)])$
and similarly
$(g \otimes g')([(f'(j'),1)]) = [(g'(f'(j')),1)]  = ((g \circ f) \otimes (g' \circ f'))([(j',1)])$. 
Since the action on edges follows trivially, $\otimes$ is indeed a bifunctor.

To see that the putative tensor product is coherent, we first note that the triangle equation turns out to be trivial, so we need only verify the pentagon equation, which we recall in Figure \ref{fig:pentagon}. The associator $\alpha_{G,G',G''} : (G \otimes G') \otimes G'' \rightarrow G \otimes (G' \otimes G'')$ is given by 
\begin{equation}
\label{eq:associator}
\alpha_{G,G',G''} :
  \begin{cases} 
    [([(j,0)],0)] & \mapsto [(j,0)] \\
    [([(j',1)],0)] & \mapsto [([(j',0)],1)] \\
    [(j'',1)] & \mapsto [([(j'',1)],1)] \\
  \end{cases}
\end{equation}
along with the implied extension to edges. The explicit form of \eqref{eq:associator} makes it clear that the associator is bijective, and hence an isomorphism.

\begin{figure}
	\centering
	\begin{tikzpicture}[->,>=stealth',shorten >=1pt]
		\node (v01) at (-4,0) {$((W \otimes X) \otimes Y) \otimes Z$};
		\node (v02) at (0,0) {$(W \otimes (X \otimes Y)) \otimes Z$};
		\node (v03) at (-4,-1.5) {$(W \otimes X) \otimes (Y \otimes Z)$};
		\node (v04) at (4,0) {$W \otimes ((X \otimes Y) \otimes Z)$};
		\node (v05) at (4,-1.5) {$W \otimes (X \otimes (Y \otimes Z))$};
		\path [->] (v01) edge node [above=1mm] {$\alpha_{W, X, Y} \otimes id_Z$} (v02);
		\path [->] (v01) edge node [left] {$\alpha_{W \otimes X, Y, Z}$} (v03);
		\path [->] (v02) edge node [above=1mm] {$\alpha_{W, X \otimes Y, Z}$} (v04);
		\path [->] (v03) edge node [above] {$\alpha_{W, X, Y \otimes Z}$} (v05);
		\path [->] (v04) edge node [right] {$id_W \otimes \alpha_{X, Y, Z}$} (v05);
	\end{tikzpicture}
	\caption{\label{fig:pentagon} The pentagon equation.}
\end{figure}

For notational convenience, let $W, X, Y, Z$ denote flow graphs with $(w, x, y, z) \in V(W) \times V(X) \times V(Y) \times V(Z)$. The three steps on the top of the pentagon are
\begin{equation}
\label{eq:pentagonThreeStep}
\begin{matrix} [([([(w,0)],0)],0)] \\ [([([(x,1)],0)],0)] \\ [([(y,1)],0)] \\ [(z,1)] \end{matrix} 
\mapsto
\begin{matrix} [([(w,0)],0)] \\ [([([(x,0)],1)],0)] \\ [([([(y,1)],1)],0)] \\ [(z,1)] \end{matrix}
\mapsto
\begin{matrix} [(w,0)] \\ [([([(x,0)],0)],1)] \\ [([([(y,1)],0)],1)] \\ [([(z,1)],1)] \end{matrix}
\mapsto
\begin{matrix} [(w,0)] \\ [([(x,0)],1)] \\ [([([(y,0)],1)],1)] \\ [([([(z,1)],1)],1)] \end{matrix}.
\end{equation}
while the two steps on the bottom of the pentagon are
\begin{equation}
\label{eq:pentagonTwoStep}
\begin{matrix} [([([(w,0)],0)],0)] \\ [([([(x,1)],0)],0)] \\ [([(y,1)],0)] \\ [(z,1)] \end{matrix} 
\mapsto
\begin{matrix} [([(w,0)],0)] \\ [([(x,1)],0)] \\ [([(y,0)],1)] \\ [([(z,1)],1)] \end{matrix}
\mapsto
\begin{matrix} [(w,0)] \\ [([(x,0)],1)] \\ [([([(y,0)],1)],1)] \\ [([([(z,1)],1)],1)] \end{matrix}.
\end{equation}
The pentagon equation follows from the equality of the rightmost parts of \eqref{eq:pentagonThreeStep} and \eqref{eq:pentagonTwoStep}, as does the theorem. \qed
\end{proof}

\section{\label{sec:Stretching}Stretching flow graphs}

By inserting new vertices and edges, we can transform many ``approximate'' flow graphs into \emph{bona fide} sub-flow graphs that can then be captured by the PST.

\begin{lemma} [``Sketch of stretch'']
\label{lemma:stretching}
Let $G$ be a flow graph. For each vertex $v \in G^\circ$, perform transformations indicated by the table below.
The cumulative result of these transformations is well-defined; repeating them has no effect. \qed
\end{lemma}

\begin{center}
\begin{tabular}{| l | c | c | c | c | c | c | c | c |}
	\hline
	$d^+(v)>1$? & $\bot$ & $\bot$ & $\bot$ & $\bot$ & $\top$ & $\top$ & $\top$ & $\top$ \\ \hline
	$d^-(v)>1$? & $\bot$ & $\bot$ & $\top$ & $\top$ & $\bot$ & $\bot$ & $\top$ & $\top$ \\ \hline
	$d^0(v)=1$? & $\bot$ & $\top$ & $\bot$ & $\top$ & $\bot$ & $\top$ & $\bot$ & $\top$ \\ \hline
	old motif 
		&
		\begin{tikzpicture}[xscale=0.4,yscale=0.5,every node/.style={inner sep=0,outer sep=0,scale=0.75},->,>=stealth',shorten >=1pt]
			\coordinate (v01) at (0,.5);
			\coordinate (v02) at (0,.25);
			\coordinate (v03) at (0,0);
			\coordinate (v04) at (1,.5);
			\coordinate (v05) at (1,.25);
			\coordinate (v06) at (1,0);
			\coordinate (v07) at (2,.5);
			\coordinate (v08) at (2,.25);
			\coordinate (v09) at (2,0);
			\foreach \from/\to in {
				v03/v06, v06/v09}
				\draw (\from) to (\to);
		\end{tikzpicture}
		&
		\begin{tikzpicture}[xscale=0.4,yscale=0.5,every node/.style={inner sep=0,outer sep=0,scale=0.75},->,>=stealth',shorten >=1pt]
			\coordinate (v01) at (0,.5);
			\coordinate (v02) at (0,.25);
			\coordinate (v03) at (0,0);
			\coordinate (v04) at (1,.5);
			\coordinate (v05) at (1,.25);
			\coordinate (v06) at (1,0);
			\coordinate (v07) at (2,.5);
			\coordinate (v08) at (2,.25);
			\coordinate (v09) at (2,0);
			\foreach \from/\to in {
				v03/v06, v06/v09}
				\draw (\from) to (\to);
			\draw (v06) to [out=60,in=120,looseness=8,min distance=20] (v06);
		\end{tikzpicture}
		&
		\begin{tikzpicture}[xscale=0.4,yscale=0.5,every node/.style={inner sep=0,outer sep=0,scale=0.75},->,>=stealth',shorten >=1pt]
			\coordinate (v01) at (0,.5);
			\coordinate (v02) at (0,.25);
			\coordinate (v03) at (0,0);
			\coordinate (v04) at (1,.5);
			\coordinate (v05) at (1,.25);
			\coordinate (v06) at (1,0);
			\coordinate (v07) at (2,.5);
			\coordinate (v08) at (2,.25);
			\coordinate (v09) at (2,0);
			\foreach \from/\to in {
				v03/v06, v06/v07, v06/v08, v06/v09}
				\draw (\from) to (\to);
		\end{tikzpicture}
		&
		\begin{tikzpicture}[xscale=0.4,yscale=0.5,every node/.style={inner sep=0,outer sep=0,scale=0.75},->,>=stealth',shorten >=1pt]
			\coordinate (v01) at (0,.5);
			\coordinate (v02) at (0,.25);
			\coordinate (v03) at (0,0);
			\coordinate (v04) at (1,.5);
			\coordinate (v05) at (1,.25);
			\coordinate (v06) at (1,0);
			\coordinate (v07) at (2,.5);
			\coordinate (v08) at (2,.25);
			\coordinate (v09) at (2,0);
			\foreach \from/\to in {
				v03/v06, v06/v07, v06/v08, v06/v09}
				\draw (\from) to (\to);
			\draw (v06) to [out=60,in=120,looseness=8,min distance=20] (v06);
		\end{tikzpicture}
		&
		\begin{tikzpicture}[xscale=0.4,yscale=0.5,every node/.style={inner sep=0,outer sep=0,scale=0.75},->,>=stealth',shorten >=1pt]
			\coordinate (v01) at (0,.5);
			\coordinate (v02) at (0,.25);
			\coordinate (v03) at (0,0);
			\coordinate (v04) at (1,.5);
			\coordinate (v05) at (1,.25);
			\coordinate (v06) at (1,0);
			\coordinate (v07) at (2,.5);
			\coordinate (v08) at (2,.25);
			\coordinate (v09) at (2,0);
			\foreach \from/\to in {
				v01/v06, v02/v06, v03/v06, v06/v09}
				\draw (\from) to (\to);
		\end{tikzpicture}
		&		
		\begin{tikzpicture}[xscale=0.4,yscale=0.5,every node/.style={inner sep=0,outer sep=0,scale=0.75},->,>=stealth',shorten >=1pt]
			\coordinate (v01) at (0,.5);
			\coordinate (v02) at (0,.25);
			\coordinate (v03) at (0,0);
			\coordinate (v04) at (1,.5);
			\coordinate (v05) at (1,.25);
			\coordinate (v06) at (1,0);
			\coordinate (v07) at (2,.5);
			\coordinate (v08) at (2,.25);
			\coordinate (v09) at (2,0);
			\foreach \from/\to in {
				v01/v06, v02/v06, v03/v06, v06/v09}
				\draw (\from) to (\to);
			\draw (v06) to [out=60,in=120,looseness=8,min distance=20] (v06);
		\end{tikzpicture}
		&
		\begin{tikzpicture}[xscale=0.4,yscale=0.5,every node/.style={inner sep=0,outer sep=0,scale=0.75},->,>=stealth',shorten >=1pt]
			\coordinate (v01) at (0,.5);
			\coordinate (v02) at (0,.25);
			\coordinate (v03) at (0,0);
			\coordinate (v04) at (1,.5);
			\coordinate (v05) at (1,.25);
			\coordinate (v06) at (1,0);
			\coordinate (v07) at (2,.5);
			\coordinate (v08) at (2,.25);
			\coordinate (v09) at (2,0);
			\foreach \from/\to in {
				v01/v06, v02/v06, v03/v06, v06/v07, v06/v08, v06/v09}
				\draw (\from) to (\to);
		\end{tikzpicture}
		&
		\begin{tikzpicture}[xscale=0.4,yscale=0.5,every node/.style={inner sep=0,outer sep=0,scale=0.75},->,>=stealth',shorten >=1pt]
			\coordinate (v01) at (0,.5);
			\coordinate (v02) at (0,.25);
			\coordinate (v03) at (0,0);
			\coordinate (v04) at (1,.5);
			\coordinate (v05) at (1,.25);
			\coordinate (v06) at (1,0);
			\coordinate (v07) at (2,.5);
			\coordinate (v08) at (2,.25);
			\coordinate (v09) at (2,0);
			\foreach \from/\to in {
				v01/v06, v02/v06, v03/v06, v06/v07, v06/v08, v06/v09}
				\draw (\from) to (\to);
			\draw (v06) to [out=60,in=120,looseness=8,min distance=20] (v06);
		\end{tikzpicture}
		 \\ \hline
	new motif 
		& 
		same
		& 
		same
		& 
		same
		& 
		\begin{tikzpicture}[xscale=0.4,yscale=0.5,every node/.style={inner sep=0,outer sep=0,scale=0.75},->,>=stealth',shorten >=1pt]
			\coordinate (v01) at (0,.5);
			\coordinate (v02) at (0,.25);
			\coordinate (v03) at (0,0);
			\coordinate (v04) at (1,.5);
			\coordinate (v05) at (1,.25);
			\coordinate (v06) at (1,0);
			\coordinate (v07) at (2,.5);
			\coordinate (v08) at (2,.25);
			\coordinate (v09) at (2,0);
			\coordinate (v10) at (3,.5);
			\coordinate (v11) at (3,.25);
			\coordinate (v12) at (3,0);
			\coordinate (v13) at (4,.5);
			\coordinate (v14) at (4,.25);
			\coordinate (v15) at (4,0);
			\foreach \from/\to in {
				v03/v06, v09/v10, v09/v11, v09/v12}
				\draw (\from) to (\to);
			\draw [dashed] (v06) to (v09);
			\draw (v06) to [out=60,in=120,looseness=8,min distance=20] (v06);
		\end{tikzpicture}
		& 
		same
		& 
		\begin{tikzpicture}[xscale=0.4,yscale=0.5,every node/.style={inner sep=0,outer sep=0,scale=0.75},->,>=stealth',shorten >=1pt]
			\coordinate (v01) at (0,.5);
			\coordinate (v02) at (0,.25);
			\coordinate (v03) at (0,0);
			\coordinate (v04) at (1,.5);
			\coordinate (v05) at (1,.25);
			\coordinate (v06) at (1,0);
			\coordinate (v07) at (2,.5);
			\coordinate (v08) at (2,.25);
			\coordinate (v09) at (2,0);
			\coordinate (v10) at (3,.5);
			\coordinate (v11) at (3,.25);
			\coordinate (v12) at (3,0);
			\coordinate (v13) at (4,.5);
			\coordinate (v14) at (4,.25);
			\coordinate (v15) at (4,0);
			\foreach \from/\to in {
				v01/v06, v02/v06, v03/v06, v09/v12}
				\draw (\from) to (\to);
			\draw (v09) to [out=60,in=120,looseness=8,min distance=20] (v09);
			\draw [dashed] (v06) to (v09);
		\end{tikzpicture}
		& 
		\begin{tikzpicture}[xscale=0.4,yscale=0.5,every node/.style={inner sep=0,outer sep=0,scale=0.75},->,>=stealth',shorten >=1pt]
			\coordinate (v01) at (0,.5);
			\coordinate (v02) at (0,.25);
			\coordinate (v03) at (0,0);
			\coordinate (v04) at (1,.5);
			\coordinate (v05) at (1,.25);
			\coordinate (v06) at (1,0);
			\coordinate (v07) at (2,.5);
			\coordinate (v08) at (2,.25);
			\coordinate (v09) at (2,0);
			\coordinate (v10) at (3,.5);
			\coordinate (v11) at (3,.25);
			\coordinate (v12) at (3,0);
			\coordinate (v13) at (4,.5);
			\coordinate (v14) at (4,.25);
			\coordinate (v15) at (4,0);
			\foreach \from/\to in {
				v01/v06, v02/v06, v03/v06, v09/v10, v09/v11, v09/v12}
				\draw (\from) to (\to);
			\draw [dashed] (v06) to (v09);
		\end{tikzpicture}
		& 
		\begin{tikzpicture}[xscale=0.4,yscale=0.5,every node/.style={inner sep=0,outer sep=0,scale=0.75},->,>=stealth',shorten >=1pt]
			\coordinate (v01) at (0,.5);
			\coordinate (v02) at (0,.25);
			\coordinate (v03) at (0,0);
			\coordinate (v04) at (1,.5);
			\coordinate (v05) at (1,.25);
			\coordinate (v06) at (1,0);
			\coordinate (v07) at (2,.5);
			\coordinate (v08) at (2,.25);
			\coordinate (v09) at (2,0);
			\coordinate (v10) at (3,.5);
			\coordinate (v11) at (3,.25);
			\coordinate (v12) at (3,0);
			\coordinate (v13) at (4,.5);
			\coordinate (v14) at (4,.25);
			\coordinate (v15) at (4,0);
			\foreach \from/\to in {
				v01/v06, v02/v06, v03/v06, v12/v13, v12/v14, v12/v15}
				\draw (\from) to (\to);
			\draw [dashed] (v06) to (v09);
			\draw [dashed] (v09) to (v12);
			\draw (v09) to [out=60,in=120,looseness=8,min distance=20] (v09);
		\end{tikzpicture}
	\\ \hline
\end{tabular}
\end{center}

Call the result of the process sketched in lemma \ref{lemma:stretching} the \emph{stretching} of $G$. This construction is similar to the ``normalization'' of two-terminal graphs (see \S \ref{sec:TTGs}).

\begin{corollary} 
There is a bijective correspondence between induced subgraphs with single sources and targets and SESE regions in a stretching. In particular, any loop corresponds to a minimal SESE region in a stretching. \qed
\end{corollary}
%
%
By considering the complete bipartite graph $K_{3,3}$, it is easy to show the following
\begin{lemma}
\label{lemma:nonplanar}
There exists a planar flow graph with a nonplanar stretching. \qed
\end{lemma}
\end{document}